\documentclass[12pt]{elsarticle}
\usepackage{supertabular}

\usepackage{epstopdf}
\usepackage{subfigure}

\usepackage[latin1]{inputenc}
\usepackage{amsfonts}

\newcommand{\sch}{Schr\"{o}dinger\ } 

\newcommand{\JCPformat}[4]{{#1} {\bf #2}, {#3} ({#4}).}

\newcommand{\Ref}[4]{\JCPformat{#1}{#2}{#3}{#4}}

\newcommand{\jcp}[3]{\Ref{J. Chem. Phys}{#1}{#2}{#3}}

\newcommand{\cpl}[3]{\Ref{Chem. Phys. Lett.}{#1}{#2}{#3}}

\newcommand{\jmathchem}[3]{\Ref{J. Math. Chem.}{#1}{#2}{#3}}
\newcommand{\jmathphys}[3]{\Ref{J. Math. Phys.}{#1}{#2}{#3}}

\begin{document}
\begin{frontmatter}

\title{An improved lower bound for the maximal length of a multivector.}
\author{Patrick Cassam-Chena\"\i}
\address{Universit\'e C\^ote d'Azur, CNRS, LJAD, UMR 7351, 06100 Nice, France.  cassam@unice.fr}

%
\date{}


\begin{abstract}
A new lower bound for the maximal length of a multivector is obtained. It is much closer to the best known upper bound than previously reported lower bound estimates. The maximal length appears to be unexpectedly large for $n$-vectors,  with $n>2$,  since the few exactly known values seem to grow linearly with vector space dimension, whereas the new lower bound has a polynomial order equal to $n-1$ like the best known upper bound. This result has implications for quantum chemistry.
\end{abstract}
\end{frontmatter}
\newpage
\section{Introduction}

The algebraic problem addressed in this article arised in the following context:  
The wave function, $\Psi$, of a quantum system made of identical Fermionic particles, such as the electrons of a molecule, or the nucleons of an atomic nucleus, is an element of the Grassmann (or Exterior) algebra, $\wedge \mathcal{H}$,  of a ``one-particle'' Hilbert space,  $\mathcal{H}$. In the non-relativistic approximation, the number of particles is a constant integer, say $n$, and the wave function actually belongs to the n$^{th}$ exterior product,   $\wedge^n \mathcal{H}$, of the Exterior algebra. In practice, in quantum physics, only approximate wave functions can be computed by solving numerically the \sch equation, as soon as $n>1$. The first step of the approach consists usually in a Galerkin-type approximation, that is to say, in reducing the infinite-dimensional, Hilbert space, $\mathcal{H}$,
 to an appropriately-chosen vector subspace, $\mathcal{H}_0$, of finite dimension, say $m$. However, such a reduction is rarely sufficient: The dimension of $\wedge^n \mathcal{H}_0$, $Dim\wedge^n \mathcal{H}_0={m \choose n}$, is still untractable, except for $m$-values  that are too small to solve accurately the time-independent \sch equation by numerical methods. So, quantum chemists define hierarchies of approximation methods that correspond to families of embedded subspaces (not necessarily vector subspaces) of $\wedge^n\mathcal{H}_0$, $V_1\subset V_2\subset\cdots\subset V_k\subset \cdots \subset\wedge^n\mathcal{H}_0$, the  $k^{th}$-order approximation corresponding to the search of an approximate solution restricted to $V_k$.\\
 
The most commonly used hierarchy is based on the concept of excitation. One starts with a zero order approximate wave function, $\Psi_0$, and a given basis sets of $\mathcal{H}_0$, $\mathcal{B}:=(\psi_1,\ldots,\psi_m)$. Then, roughly speaking, i.e. omitting subtilities related to particle spin and other symmetries, $V_k$ is the vector space spanned by all $n$-vectors built by substituting in the expansion of $\Psi_0=\sum\limits_Ia_I\Psi_I$ in the induced basis set, $\wedge^n\mathcal{B}:=(\Psi_I)_{I=i_1<\cdots <i_n}:=(\psi_{i_1}\wedge\cdots\wedge\psi_{i_n})_{i_1<\cdots <i_n}$, at most $k$ elements of $\mathcal{B}$ by $k$ other such basis elements. $V_k$ is the space of the so-called ``$l$-excited configurations (from  $\Psi_0$)'' with $l\leq k$ \cite{Szabo89}.

Another common hierarchy is based on the rank of an $n$-vector, $r(\Psi)$ defined as the least integer $l$ such that $\exists \mathcal{F}\subseteq \mathcal{H}_0,\ Dim\,\mathcal{F}=l$ and $\Psi\in\wedge^n\mathcal{F}$. In this case, $V_k=\{\Psi\in\wedge^n\mathcal{H}_0, r(\Psi)\leq k\}$, that is the space of  $n$-vectors of rank less than or equal to $k$. The method which consists in finding the best wave function in the variational sense, that is to say, by minimizing the energy functional of the system, subject to the constraint that the wave function belongs to $V_k$, is called the ``complete active space self-consistent field'' method (with active space of dimension $k$) \cite{Cassam92-lma,Cassam94-jmc}. 

A third possible hierarchy is based on the so-called ``depth'' of an $n$-vector, $d(\Psi)$ defined as the largest integer $l$ such that 
$\exists \Psi_1\in \wedge^{n_1} \mathcal{H}_0,\ldots,\exists\Psi_l\in \wedge^{n_l} \mathcal{H}_0$,
with $n_1,\ldots ,n_l$ non negative integers, and $\Psi=\Psi_1\wedge\cdots\wedge\Psi_l$ \cite{Cassam03-jmp}. If one sets, $V_k=\{\Psi\in\wedge^n\mathcal{H}_0, d(\Psi)\geq k\}$, a hierarchy of approximations corresponds to the embedding,  $V_n\subset \cdots\subset V_1=\wedge^n\mathcal{H}_0$. The variational solution of maximal depth, i.e. $\Psi\in V_n$, is actually the unrestricted Hartree-Fock  \cite{Berthier54,Pople54} wave function, where all its factors $\Psi_i\in\mathcal{H}_0\equiv\wedge^1\mathcal{H}_0$. The case where all  factors $\Psi_i\in\wedge^2\mathcal{H}_0$ for $n$ even, has been explored in \cite{Cassam06-jcp,Cassam07-cpl,Cassam10-cpl}.

These hierarchies can give rise to interesting algebraic problems. However, the problem we are dealing with in this paper is related to yet another hierchachy based on the length of an $n$-vector. The length of an $n$-vector, $l(\Psi)$ is the least positive integer $l$, such that,
$\exists \Psi_1\in \wedge^{n} \mathcal{H}_0, \ldots, \exists\Psi_l\in\wedge^{n}\mathcal{H}_0$, with $\Psi_1,...,\Psi_l$ decomposable (also called ``pure'' or ``simple'') $n$-vectors, and
$\Psi=\Psi_1+\cdots+\Psi_l$. Note that generalizations of this concept and their connection to quantum physics can be found in \cite{Cassam04-pla}.
The length-based hierarchy has been proposed in the early 90's \cite{Cassam92,Cassam94-jmc}. 
In this case, $V_k=\{\Psi\in\wedge^n\mathcal{H}_0, l(\Psi)\leq k\}$, that is the space of $n$-vectors of length less than or equal to $k$. The variational solution of minimal length, i.e. $\Psi\in V_1$, is again the unrestricted Hartree-Fock wave function. 
More recently, Beylkin and coworkers have worked out  techniques to overcome the numerical difficulties (mostly related to non-orthogonality of the elements of $\mathcal{H}_0$ making up the decomposable $\Psi_k$ terms) in optimizing wave functions of increasing length \cite{Beylkin08}.

However,  to be really effective, such a method must be able to give accurate approximate wave functions of reasonably small length. There is \textit{a priori} no reason for the best approximate wave function $\Psi\in\wedge^n\mathcal{H}_0$ to have a length less than the maximal length achievable by an $n$-vector. Hence, the interest of determining such a maximal length, which is still an open problem. This topics has been reviewed by MacDougall \cite{MacDougall82}. The maximal length depends not only upon $n$ and $m$, but also upon the field $\mathbb{K}\equiv\wedge^0\mathcal{H}_0$ \cite{Busemann73}, so we will use the notation  $N( \mathbb{K},m,n)$ for it. In quantum physics , the relevant field is usually taken to be the field of complex numbers $\mathbb{C}$. However, in many occasions, one uses the field of real numbers $\mathbb{R}$, and, in some instances, the field of quaternions, $\mathbb{H}$, can also be considered. When a result is independent of the field, we will write simply, $N(m,n)$. 

The purpose of this paper is to give a new lower bound to $N( \mathbb{C},m,n)$, which improves drastically those reported previously. This is achieved in the next section, before we conclude in the final section.

\section{Lower bound to the maximal length}

What is known about the maximal length $N( \mathbb{K},m,n)$ is essentially summarized in \cite{MacDougall82}.
First, it is useful to note that $N(m,n)=N(m,m-n)$ by what would be called in quantum physics the ``particle-hole duality''
between $\wedge^n\mathcal{H}_0$ and $\wedge^{m-n}\mathcal{H}_0$. So, we can limit the study to $N(m,n),\ n\leq \lfloor\frac{m}{2}\rfloor$, denoting by $\lfloor k\rfloor$ the floor value of $k$. Exact values are known in very few cases:
For the limiting case $n=1$, obviously $N(m,1)=1$. For $n=2$, one obtains easily by making use of Schmidt decomposition that $N(m,2)=\lfloor\frac{m}{2}\rfloor$.  For $n=3$, Glassco showed that $N(\mathbb{K},6,3)=3$ for any $\mathbb{K}$ of characteristic $0$ \cite{MacDougall82}; Westwick showed that $N(\mathbb{C},7,3)=4$ and that $N(\mathbb{C},8,3)=5$. This is all, to our knowledge.\\

Exact values being difficult to obtain, one has tried to find bounds.
In the general case, i.e. for arbitrary values of $n$ and $m$, the best known lower  bound to $N(m,n)$  is  $\lfloor \frac{m-n+2}{2} \rfloor $, \cite{Busemann73}, whereas the best upper bound is a polynomial  of the order of $\frac{m^{n-1}}{2n!}$ \cite{MacDougall82}. So, the best known lower and upper  bounds are  extremely far apart. Furthermore, for a fixed $m$, the lower bound is exact for $n=2$, then it decreases with $n$, whereas $N(m,n)$ is expected to increase with $n$ as long as $n\leq \lfloor\frac{m}{2}\rfloor$. We will now show how to find a lower bound that remedies to both of these discrepancies for $n>2$.\\

Let us assumed that $\mathbb{K}=\mathbb{C}$, so that the field is algebrically closed.
The length of any  $\Psi\in\wedge^n\mathcal{H}_0$  being less than or equal to $N(\mathbb{C},m,n)$, (by definition of $N(\mathbb{C},m,n)$,)  we can write, $\forall \Psi\in\wedge^n\mathcal{H}_0, \exists (\phi^1_{i_1}, \ldots ,\phi^1_{i_n}, \ldots , \phi^{N(\mathbb{C},m,n)}_{i_1},\ldots ,\phi^{N(\mathbb{C},m,n)}_{i_n})\in\mathcal{H}_0^{n\cdot N(\mathbb{C},m,n)}$, such that,
\begin{eqnarray}
\Psi=\sum\limits_{j=1}^{N(\mathbb{C},m,n)} \phi^j_{i_1}\wedge \cdots \wedge\phi^j_{i_n}.
\label{expansion-length}
\end{eqnarray} 
This expresses that $N(\mathbb{C},m,n)$ is the smallest integer, $k$, such that the union for $l<k$ of the $l$-secant varieties, $S_l(G(n,m))$,  of the Grassmannian, $G(n,m)$, contains the whole of $\wedge^n\mathcal{H}_0$.

Recall that the $l$-secant variety,  of a non-degenerate projective variety $X$, $S_l(X)$, is usually defined to be the closure of the union of linear spans of all $(l+1)$-uples of independent points lying on $X$ \cite{Harris1995,Ciliberto2008}. Note however, that in a recent reference concerned with secants of Grassmannian \cite{Boralevi2013}, the $l$-secant is defined to be the closure of the union of linear spans of all $(l)$-uples of independent points, so that it is actually what we would denote $S_{l-1}(X)$.

The projective dimension of the Grassmannian $G(n,m)$ is known to be $n(m-n)$ and the projective dimension of $\wedge^n\mathcal{H}_0$ in which the $G(n,m)$ is embedded is ${m \choose n} -1$. 
The dimension of the $(l-1)$-secant variety is bounded by \cite{Ciliberto2008},
\begin{eqnarray}
Dim\,S_{l-1}(G(n,m)) \leq Min\{ l\cdot \left(n\cdot(m-n)+1\right)-1,{m \choose n} -1\}.
\label{secant-dim}
\end{eqnarray} 
The $(l-1)$-secant variety is said ``to have the expected dimension'' when the equality holds, otherwise it is said ``defective'', its defect being the difference between the upper bound in Ineq.~(\ref{secant-dim}) and its actual dimension.
For $l=N(\mathbb{C},m,n)$, Eq.~(\ref{expansion-length}) implies that this dimension cannot be less than the projective dimension of $\wedge^n\mathcal{H}_0$,
so, ${m \choose n} -1 \leq N(\mathbb{C},m,n)\cdot\left(n\cdot(m-n)+1\right)-1$, hence
\begin{eqnarray}
\frac{{m \choose n}}{n\cdot(m-n)+1}  \leq N(\mathbb{C},m,n).
\end{eqnarray}

\begin{figure}[h]
\begin{center}
\includegraphics[scale=0.6]{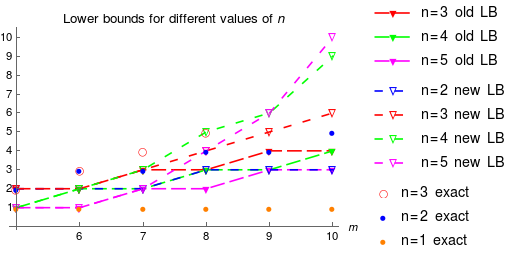}
\caption{Comparison of the new lower bounds (empty triangles) with respect to the old ones (filled triangles) for different values of $m$ and $n$. Known exact values of $N(\mathbb{C},m,n)$ are also plotted with empty ($n=3$) or filled ($n<3$) circles. Different colors correspond to different $n$-values.} 
 \label{LB-fig}
\end{center}
\end{figure}

Note that this new lower bound is \textit{a fortiori} a lower bound for $N(\mathbb{R},m,n)$. Not surprisingly, it is not very tight for $n=2$, because in this case $S_{l-1}(G(n,m))$ is almost always defective \cite{Boralevi2013}. However, it is much tighter than $\lfloor \frac{m-n+2}{2} \rfloor $ for $n>2$ (see Fig.~\ref{LB-fig}) and increases with $n$  as long as $n\leq \lfloor\frac{m}{2}\rfloor$, as expected.
It is manifestly invariant by  duality between $\wedge^n\mathcal{H}_0$ and $\wedge^{m-n}\mathcal{H}_0$.
For $m$ very large compared to $n$, it is in the order $\frac{m^{n-1}}{n\cdot n!}$ so
within a ratio $\frac{2}{n}$ with respect to the order of the known upper bound mentionned above.

\newpage
\section{Conclusion}

The lower bound to the length of an $n$-vector obtained in this work is much higher than was anticipated \cite{MacDougall82}. 
It turns out that the few known values of $N(\mathbb{C},m,n)$ were deceptively low. In fact, the known upper bound for large $m$,
was much better than previously thought \cite{MacDougall82}, as it is quite close to our new lower bound.
\textit{A priori}, this is not a good news for length-based approximations in quantum chemistry, because there is no reason why
the best approximate wave function in  $\wedge^n\mathcal{H}_0$  for a molecular system should have a length less than $N(\mathbb{C},m,n)$.
Based on the known values of $N(\mathbb{C},m,n)$, one could have hoped that the computational cost of a length-based hierarchy of approximations would scale linearly with $m$, we deduce from our new lower bound that it should rather scale as $m^{n-1}$. 
However, it can still be hoped that for physical reasons, only a few terms in its expansion, Eq.~(\ref{expansion-length}), would contribute significantly to wave function, once normalized by its $L^2$-norm. 

The problem of the length of an $n$-vector does not seem to have attracted much attention in the mathematical community since the 1970's. We hope that the present paper will revive interest in this topics and in similar mathematical problems related to the depth of an $n$-vector. Progress in this field are potentially important for quantum physics, in particular for quantum chemistry approximation methods, or quantum entanglement related technologies. 

\section*{Acknowledgements}
B. Mourrain and A. Galligo, are ackowledged for discussions that have helped the author in putting
his arguments in the present form.  In particular, B. Mourrain has provided precise algebraic geometry
hints and references.

\end{document}